\documentclass{article}
\usepackage{amsmath, amssymb, amsthm}

\usepackage[dvips]{graphicx}
\usepackage[dvips]{epsfig}
\usepackage{psfrag}
\usepackage{txfonts}

\theoremstyle{definition}

\title{ K-theory of some C*-algebras and buildings.}

\author{\sc Alina Vdovina}

\date{}

\begin{document}

\maketitle

\abstract{ We compute an exact formula for the order of the class 
of the identity
in the $K_0$ group of an infinite class of two-dimensional
Kuntz-Crieger algebras.}

\section*{Introduction}

The class of the identity [$\mathbf{1}$] in $K_0$
of different classes of crossed product
 $C^*$-algebras was broadly investigated
in the literature, see \cite{Co1}, \cite{N}, \cite{HN}, \cite{AD},
\cite{RS2}, \cite{M}.

We will concentrate on a case associated to two-dimensional
Euclidean buildings.
Let group $G$ acts simply transitively on the vertices 
of a $\Tilde{A}_2$ building
$\Delta$. Then there is an induced action on the boundary $\Omega$ of $\Delta$,
the crossed product algebra $C(\Omega) \rtimes \Gamma$ depends only on $G$
and is classified by its $K$-groups together 
with the class [$\mathbf{1}$] in $K_0$
of the identity element $\mathbf{1}$ of $C(\Omega) \rtimes \Gamma$.
It is interesting therefore to identify this class.

We will consider special class of $\Tilde{A}_2$ 
groups $\Gamma_{\mathcal{T}_0}$ described in \cite{Cart},
which embed as arithmetic subgroups
of $PGL(3,\mathbf{F}_q(X))$. For this class of groups we prove 
the following result,
which was conjectured for all $\Tilde{A}_2$ groups in \cite{RS2}.

\bigbreak

{\bf Theorem.} The order of the class [$\mathbf{1}$] of the identity element
$\mathbf{1}$ of

 $C(\Omega) \rtimes \Gamma$ in $K_0(C(\Omega) \rtimes \Gamma)$ is $q-1$,
where $\Gamma$ is a  $\Gamma_{\mathcal{T}_0}$ group and $q \nequiv 1 (mod~3)$.

\bigbreak

\section*{Polygonal presentation and construction of polyhedra.}

A {\em polyhedron} is a
two-dimensional complex which is obtained from several oriented
$p$-gons by identification of corresponding sides.
 Let's take  take a sphere  of
a small radius at a point of the polyhedron. 
The intersection of the sphere with
the polyhedron is a graph, which is called the {\em link} at this
point. We consider a case, when 
 all sides of a polyhedron are regular euclidean triangles
and links at all vertices are incidence graphs of finite projective planes.
The universal covering of such a polyhedron is an euclidean
 $\Tilde{A}_2$ building \cite{BB}, \cite{Ba} and
 with the metric
introduced in \cite[p.~165]{BBr} it is a complete metric space
of non-positive curvature in the sense of Alexandrov and
Busemann. 
 It follows from \cite{BS}, that
the fundamental groups of our polyhedra satisfy
the property (T) of Kazhdan. (Another relevant reference
is \cite{Z}.)

\medbreak

\noindent
{\bf Definition.} Let $\mathcal{P}$ be a tesselation of the
Euclidean plane by  regular triangles,
with angles $\pi/3$ in each vertex.
A {\em Euclidean} $\Tilde{A}_2$ {\em building} is a polygonal complex $X$,
which can be expressed as the union of subcomplexes called apartments
such that:

\begin{itemize}
\item[1.] Every apartment is isomorphic to $\mathcal{P}$.
\item[2.] For any two polygons of $X$, there is an apartment
containing both of them.
\item[3.] For any two apartments $A_1, A_2 \in X$ containing
the same polygon, there exists an isomorphism $ A_1 \to A_2$
fixing $A_1 \cap A_2$.
\end{itemize}

\medbreak

Recall  that
a {\em generalized m-gon} is a connected, bipartite graph of
diameter $m$ and girth $2m$, in which each vertex lies on at least
two edges. A graph is {\em bipartite} if its set of vertices can
be partitioned into two disjoint subsets such that no two vertices
in the same subset lie on a common edge. The vertices of
 one subset we will call black vertices and the
vertices of the other subset the white ones. The {\em diameter} is
the maximum distance between two vertices and the {\em girth}  is
the length of a shortest circuit. Incidence graphs of finite
projective planes are exactly generalized triangles.

We recall a definition of polygonal presentation
introduced in \cite{V}.

\medbreak

\noindent {\bf Definition.} Suppose we have $n$ disjoint connected
bipartite graphs\linebreak
 $G_1, G_2, \ldots G_n$.
Let $P_i$ and $Q_i$ be the sets of black and white vertices
respectively in $G_i$, $i=1,\dots,n$; let $P=\bigcup P_i$,
$Q=\bigcup Q_i$, $P_i \cap P_j = \emptyset$,
 $Q_i \cap Q_j = \emptyset$
for $i \neq j$ and
let $\lambda$ be a  bijection $\lambda: P\to Q$.

A set $\mathcal{K}$ of $k$-tuples $(x_1,x_2, \ldots, x_k)$, $x_i \in P$,
will be called a {\em polygonal presentation} over $P$ compatible
with $\lambda$ if

\begin{itemize}

\item[(1)] $(x_1,x_2,x_3, \ldots ,x_k) \in \mathcal{K}$ implies that
   $(x_2,x_3,\ldots,x_k,x_1) \in \mathcal{K}$;

\item[(2)] given $x_1,x_2 \in P$, then $(x_1,x_2,x_3, \ldots,x_k) \in \mathcal{K}$
for some $x_3,\ldots,x_k$ if and only if $x_2$ and $\lambda(x_1)$
are incident in some $G_i$;

\item[(3)] given $x_1,x_2 \in P$, then  $(x_1,x_2,x_3, \ldots ,x_k) \in \mathcal{K}$
    for at most one $x_3 \in P$.

\end{itemize}

If there exists such $\mathcal{K}$, we will call $\lambda$ a {\em
  basic bijection}.

\medskip

We can associate  a polyhedron $K$ on $n$ vertices with
each polygonal presentation $\mathcal{K}$ as follows:
for every cyclic $k$-tuple $(x_1,x_2,x_3,\ldots,x_k)$ 
we take an oriented $k$-gon on the boundary of which
the word $x_1 x_2 x_3\ldots x_k$ is written. To obtain
the polyhedron we identify the corresponding sides of our
polygons, respecting orientation.

\medskip

\noindent {\bf Lemma \cite{V}} A polyhedron $K$ which
corresponds to a polygonal presentation $\mathcal{K}$ has
  graphs $G_1, G_2, \ldots, G_n$ as the links.

\medskip

\section{Balanced polygonal presentation}

We will use a particular case of polygonal presentation,
so-called triangle presentation, described in \cite{Cart}.
We repeat now the construction from \cite{Cart} for completness.

Consider the Desarguesian projective plane $(P,L)=PG(2,q)$
of prime power order $q$, in which the points and lines
are 1- and 2-dimensional subspaces, respectively,
of a 3-dimensional vector space $V$ over $\mathbf{F}_q$,
with incidence being inclusion. We may take $V=\mathbf{F}_{q^3}$.
Consider a regular quadratic form on $F_{q^3}$ $(x_0,y_0) \to Tr(x_0,y_0)$,
where $Tr$ is the trace of the field extension $F_{q^3}/F_q$.
For $x \in P$, set
$$ \lambda_0=\{y \in P: Tr(x,y)=0\}.$$

This defines a point-line correspondence $\lambda_0: P \to L$.
The following set of triples is a triangle presentation 
$\mathcal{T}_0$ compatible with $\lambda_0$:
$$\mathcal{T}_0=\{(x, x\xi, x\xi^{q+1})|x, \xi \in P, Tr(\xi)=0\}.$$

\noindent
It is convinient to denote elements of $P$ with letters of
an alphabet $\mathcal{X}=\{x_1, \dots, x_{q^2+q+1}\}$.

We describe now a new polygonal presentation $\mathcal{T}_1$.
 Take an alphabet 
$\mathcal{Y}$=$\mathcal{A}, \mathcal{B}, \mathcal{C}$,
were every subalphabet $\mathcal{A}, \mathcal{B}, \mathcal{C}$
contains $q^2+q+1$ elements. 
$\mathcal{A}=\{a_i\}$, $\mathcal{B}=\{b_i\}$, 
$\mathcal{C}=\{c_i\}$, $i=1, \dots, q^2+q+1$.

Define $\mathcal{T}_1$ as the following set of triples:
$$\{(a_k,b_l,c_m),\, (b_k, c_l, a_m),\,
(c_k,a_l, b_m) \in \mathcal{T}_1 \iff (x_k, x_l, x_m) \in \mathcal{T}_0\}$$

Define bijections $\lambda_1$,
$\lambda_2$,$\lambda_3$
in the following way
$\lambda_1(x_i)=a_i, \,\lambda_2(x_i)=b_i, \,\lambda_3(x_i)=c_i$.

\bigbreak

{\bf Lemma 1.} There exists a subset of  $\mathcal{T}_1$, such
that every element of $\mathcal{Y}$ occurs exactly once.
We will call such a subset $\mathcal{S}$ {\em basic} subset.
A polygonal presentation containing a basic subset will be
called {\em balanced} presentation.

\bigbreak

{\bf Proof.} 
Let's consider such an element $x \in P$, that is a generator
of $P$ as a  cyclic group. Fix $\xi \in P$
 such that $Tr(\xi)=0$.  Now, consider
the following set of triples $(a_i,b_i,c_i), i=1,...,g^2+q+1$,
where $a_i=\lambda_1(x^i)$, $b_i=\lambda_2(x^i\xi)$,
$c_i=\lambda_3(x^i \xi^{q+1})$.

\section{Subshift of a balanced polygonal presentation}

Let  $\mathcal{T}$ be a polygonal presentation with $n=3$, $k=3$,
where all there graphs $G_1$, $G_2$ and $G_3$ are incidence
graphs of finite projective planes of order $q$. 
The polyhedron, which corresponds
to $\mathcal{T}$, has triangular faces and three vertices.
We will consider polyhedra such that all three vertices of each 
triangle have different graphs as links.
In this case  
we can give a Euclidean metric to every face. In this metric
 all sides of the triangles
are geodesics of the same length. The universal covering of the polyhedron
 is an Euclidean building $\Delta$, see \cite{BB}, \cite{Ba}.
 Each element of $\mathcal{T}$ may be identified with  
an oriented basepointed triangle in $\Delta$. We now construct
a 2-dimensional shift system associated with  $\mathcal{T}$.
The transition matrices $M_1, M_2$
in the way, defined as in \cite{RS2},p.828: if $\alpha=(x_1,x_2,x_3),
 \beta=(y_1, y_2, y_3) \in \mathcal{T}$
say that $M_1( \beta, \alpha)=1$ if and only if there exists
  $\psi=(x_3,z,y_1)$ and $M_1( \beta, \alpha)=0$ otherwise (Figure \ref{fig-1}).
 In a similar way,
$M_2(\gamma, \alpha)=1$ for $\alpha=(x_1,x_2,x_3), \gamma=(y_1,y_2,y_3)$
if and only if there exists $\psi=(x_2,y_1,z)$ and
$M_2(\gamma, \alpha)=0$
otherwise.

\begin{figure}
\setlength{\unitlength}{1mm}
\begin{picture}(70,50)
\put(20,5){\vector(1,2){5}}
\put(25,15){\line(1,2){5}}
\put(30,25){\vector(1,-2){5}}
\put(35,15){\line(1,-2){5}}
\put(20,5){\line(1,0){10}}
\put(40,5){\vector(-1,0){10}}
\put(29,12){$\alpha$}
\put(28,2){$x_1$}
\put(21,15){$x_2$}
\put(36,14){$x_3$}
\put(30,25){\vector(1,2){5}}
\put(35,35){\line(1,2){5}}
\put(40,45){\vector(1,-2){5}}
\put(45,35){\line(1,-2){5}}
\put(30,25){\line(1,0){10}}
\put(50,25){\vector(-1,0){10}}
\put(39,32){$\beta$}
\put(38,22){$y_1$}
\put(31,35){$y_2$}
\put(46,34){$y_3$}
\put(40,5){\vector(1,2){5}}
\put(45,15){\line(1,2){5}}
\put(45,12){$z$}
\put(95,5){\vector(-1,0){10}}
\put(85,5){\line(-1,0){10}}
\put(75,5){\vector(1,2){5}}
\put(80,15){\line(1,2){5}}
\put(85,25){\vector(1,-2){5}}
\put(90,15){\line(1,-2){5}}
\put(84,12){$\alpha$}
\put(83,2){$x_1$}
\put(76,15){$x_2$}
\put(91,14){$x_3$}
\put(85,25){\vector(-1,0){10}}
\put(75,25){\line(-1,0){10}}
\put(65,25){\vector(1,2){5}}
\put(70,35){\line(1,2){5}}
\put(75,45){\vector(1,-2){5}}
\put(80,35){\line(1,-2){5}}
\put(74,32){$\gamma$}
\put(73,22){$y_1$}
\put(66,35){$y_2$}
\put(81,34){$y_3$}
\put(65,25){\vector(1,-2){5}}
\put(70,15){\line(1,-2){5}}
\put(68,13){$z$}
\end{picture}
\label{fig-1}
\caption{}
\end{figure}

 The matrices  $M_1, M_2$ of order
$3(q+1)(q^2+q+1) \times 3(q+1)(q^2+q+1)$ are nonzero
$\{0,1\}$ matrices.
We will use $\mathcal{T}$ as an alphabet and $M_1, M_2$
as transition matrices to build up 2-dimensional words as in
\cite{RS1}. Let $[m,n]$ denote $\{m, m+1,...,n\}$, where $m \leq n$
are integers. If $m,n \in \mathbb{Z}^2$, say that $m \leq n$
if  $m_j \leq n_j$ for j=1,2, and when $m \leq n$ let
$[m,n]=[m_1,n_1] \times [m_2,n_2]$. In $ \mathbb{Z}^2$, let 0
denote the zero vector and let $e_j$ denote the $j$-th standard 
unit basis vector. If $m \in \mathbb{Z}_+^2=\{m \in \mathbb{Z}^2;
m \geq 0\}$, let

\noindent $W_m=\{w:[0,m] \to \mathcal{T}; M_j(w(l+e_j),w(l)=1$ where
$l,l+e_j \in [0,m]\}$
and call the elements of $W_m$ words.

 In order to apply the theory
from \cite{RS1} we need the matrices $M_1,M_2$ to satisfy the following
conditions:

\bigbreak

(H0) Each $M_i$ is a nonzero $\{0,1\}$-matrix.

(H1a) $M_1M_2=M_2M_1$.

(H1b) $M_1M_2$ is a $\{0,1\}$-matrix.

(H2) The directed graph with vertices $\alpha \in \mathcal{T}$
and directed edges $(\alpha,\beta)$
whenever $M_i(\alpha, \beta)=1$ for some $i$ is irreducible.

(H3) For any nonzero $p \in \mathbb{Z}^2$, there exists
a word $w \in W$ which is not $p-periodic$, i.e., there exists $l$
so that $w(l)$ and $w(l+p)$ are both defined but not equal.

\bigbreak

In \cite{RS1} some $C^*$-algebra is defined by partial isometries of
 the system of words $W_m$, where $m \in \mathbb{Z}_+^2$.
It is proved there, that if the matrices  
$M_1,M_2$  satisfy the conditions
(H0),(H1a,b),(H2),(H3), then this algebra is simple,
purely infinite and nuclear.

\medbreak

Now we prove the conditions (H0), (H1a,b), (H2), (H3)
for our two-dimensional shift.
By definition of matrices $M_1,M_2$ they are nonzero
$\{0,1\}$ matrices, so (H0) holds. 
If we have $\alpha$, $\beta$, $\psi$, such that  $M_1(\alpha, \beta)=1$,
$M_2(\beta, \psi)=1$, then $\gamma$ such 
that $M_2(\alpha, \gamma)=1$, $M_1(\gamma, \psi)=1$, is uniquely
defined because of properties
of finite projective planes. Conditions (H1a,b) follow.
To prove (H2) we need to color sides of triangles in three
different colors. This is possible since there are three vertices
in the polyhedron with different graphs as links.
So, all triangles from $\mathcal{T}$ have one of three possible
colorings.
We need to show, that for any $\alpha, \beta \in \mathcal{T}$
we can choose $r>0$ such that $M_j^r(\alpha, \beta)>0$, where $j=1,2$.
Geometrically it means that any  $\alpha, \beta \in \mathcal{T}$
can be realized so that $\beta$ lies in some sector with base $\alpha$
(for more detailes see \cite{RS1}). Without loss of generality
we can assume, that $j=1$. We will say, that $\beta \in \mathcal{T}$
is reachable from $\alpha \in \mathcal{T}$ in $r$ steps, if there is
 $r>0$ such that $M_1^r(\alpha, \beta)>0$. It is easy to see, that
every triangle is reachable from some triangle of other color
in one or two steps. So, to prove (H2) we need to show, that
any triagnle is reachible from another one of the same color.
Now we can use the proof of the Theorem 1.3 from \cite{RS1},
since at each step of this proof it is only used, that
the link at each vertex of the building is an incidence graph
of a finite projective plane, which is true in our case too.
The proof of (H3) is identical to the proof of (H3)
in the case of the subshift considered in \cite{RS3}.

\bigbreak

Now, as a set of triangles we consider all elements of  $\mathcal{T}_1$,
every cyclic word $(a_i, b_j, c_k)
 \in \mathcal{T}_1$ brings three basepointed triangles (Figure \ref{fig-2}).

\begin{figure}
\setlength{\unitlength}{1mm}
\begin{picture}(70,30)
\put(40,5){\vector(-1,0){10}}
\put(30,5){\line(-1,0){10}}
\put(20,5){\vector(1,2){5}}
\put(25,15){\line(1,2){5}}
\put(30,25){\vector(1,-2){5}}
\put(35,15){\line(1,-2){5}}
\put(28,2){$a_i$}
\put(21,15){$b_j$}
\put(36,14){$c_k$}
\put(70,5){\vector(-1,0){10}}
\put(60,5){\line(-1,0){10}}
\put(50,5){\vector(1,2){5}}
\put(55,15){\line(1,2){5}}
\put(60,25){\vector(1,-2){5}}
\put(65,15){\line(1,-2){5}}
\put(58,2){$b_j$}
\put(51,15){$c_k$}
\put(66,14){$a_i$}
\put(100,5){\vector(-1,0){10}}
\put(90,5){\line(-1,0){10}}
\put(80,5){\vector(1,2){5}}
\put(85,15){\line(1,2){5}}
\put(90,25){\vector(1,-2){5}}
\put(95,15){\line(1,-2){5}}
\put(88,2){$c_k$}
\put(81,15){$a_i$}
\put(96,14){$b_j$}
\end{picture}
\label{fig-2}
\caption{}
\end{figure}

{\bf Lemma 2.} The set $M(\mathcal{S}^a)$ consists of $q-1$ copies of
every element of  $\mathcal{T}_1^b$ and one copy of  $\mathcal{S}^b$.

\medbreak

Denote $\mathcal{S}^a$ tiles of $\mathcal{S}$ starting with $a_i$,
$i=1,...,g^2+q+1$
and analyse, which elements of  $\mathcal{T}_1$, and how
many of them can be obtained by one left shift from $\mathcal{S}^a$.
In general, from each tile in $\Tilde{A}_2$ case, one
can obtain $q^2$ tiles by one left(right) shift
as a consequence of properties of finite projective planes
(see [RS] for details). Now, each tile $\gamma \in \mathcal{T}_1^b$
can be obtained from some tile $\alpha \in \mathcal{S}^a$
by $q^2$ times. So, the total number of tiles which can be obtained
from $\mathcal{S}^a$ is $q^2(q^2+q+1)$. Since every $c_i$ appears
exactly once, every $\gamma \in \mathcal{T}_1^b$ will
appear in $M(\mathcal{S}^a)$  exactly $q$ times if
 $\gamma \in \mathcal{S}_1^b$ and $q-1$ otherwise.
So, the set $M(\mathcal{S}^a)$ consists of $q-1$ copies of
every element of  $\mathcal{T}_1^b$ and one copy of  $\mathcal{S}^b$.

\medbreak 

Two subsequent copies of this lemma one gets by substitution
$a$ by $b$, $b$ by $c$ and $a$ by $c$, $b$ by $a$.

\section{The class of the identity in $K$-theory.}


{\bf Proof of the Theorem.}

It was shown in \cite{RS2}, that the $K$-theory of the $C^*$-algebra
can be found as the abelian group with generators,
which are the elements of the alphabet $\mathcal{T}_1^b$
with the following relations:
 $$t=\sum_{s \in \mathcal{T}_1} M_1(s,t)$$

It follows from \cite{CMS}, \cite{RS2}, that the identity 
function in $C(\Omega)$
can be expressed as the sum of all tiles 
of the alphabet  $\mathcal{T}_1$, $$\sum_{t \in \mathcal{T}_1} t \,\,\, .$$

So, we will use the system of relations
$$t=\sum_{s \in \mathcal{T}_1} M_1(s,t)$$
to express  $$\sum_{t \in \mathcal{T}_1} t \,\,\, . $$

It follows from Lemma 2, that 
 $$\sum_{t \in \mathcal{S}^a} t = (q-1)\sum_{t \in \mathcal{T}_1^b} t +
  \sum_{t \in \mathcal{S}^b} t$$

$$\sum_{t \in \mathcal{S}^b} t = (q-1)\sum_{t \in \mathcal{T}_1^c} t +
  \sum_{t \in \mathcal{S}^c} t$$

$$\sum_{t \in \mathcal{S}^c} t = (q-1)\sum_{t \in \mathcal{T}_1^a} t +
  \sum_{t \in \mathcal{S}^a} t \,\,\, .$$

\noindent
By addition of these three equalities we get

$$(q-1)\sum_{t \in \mathcal{T}_1} t = 0 \, ,$$

\noindent
so $(q-1)\mathbf{1}=1$. But it was shown in
\cite{RS2}, that the order of $\mathbf{1}$ is at least $q-1$ 
in the case when $q \nequiv 1 (mod~3)$,
what completes the proof.

\medbreak
\noindent
{\em Acknowledgement}.

\noindent
I would like to thank G.~Robertson for useful discussions
and comments.

\end{document}